\documentclass[reqno, a4paper,12pt]{amsart}

\usepackage[capposition=bottom]{floatrow}
\usepackage{euscript,eufrak,verbatim}
\usepackage{graphicx}
\usepackage[usenames]{color}
\usepackage[colorlinks,linkcolor=red,anchorcolor=blue,citecolor=blue]{hyperref}
\usepackage{amsmath, mathtools}
\usepackage{mathtools}
\usepackage{amsthm}
\usepackage[all]{xy}
\usepackage{amssymb} 
\usepackage{bm}
\usepackage{booktabs,siunitx}

\usepackage[all]{xy}

\usepackage{mathrsfs}
\usepackage{amscd}

\usepackage{enumitem}
\usepackage[numbers]{natbib}

\usepackage{mathptmx}
\usepackage[T1]{fontenc}

\usepackage{lipsum}
\usepackage{accents}
\usepackage{titlesec}

\makeatletter
\numberwithin{equation}{section}

\theoremstyle{cuplain}
\newtheorem{main theorem}{Main Theorem}
\newtheorem{theorem}{Theorem}[section]
\newtheorem{lemma}[theorem]{Lemma}
\newtheorem{conjecture}[theorem]{Conjecture}

\newtheorem{proposition}[theorem]{Proposition}

\newtheorem*{theorem*}{``Theorem''}
\theoremstyle{definition}
\newtheorem{definition}[theorem]{Definition}
\newtheorem{remark}[theorem]{Remark}
\newtheorem{example}[theorem]{Example}

\newtheorem*{example*}{Example}
\newtheorem*{remark*}{Remark}

\newtheoremstyle{break}
  {\topsep}{\topsep}%
  {\itshape}{}%
  {\bfseries}{}%
  {\newline}{}%
\theoremstyle{break}

\newtheoremstyle{break}
  {\topsep}{\topsep}%
  {\normalshape}{}%
  {\bfseries}{}%
  {\newline}{}%

\newtheorem{breakremark}[theorem]{Remark}

\newenvironment{retheorem}[1]{%
  \par\medskip\noindent
  \textbf{Theorem \ref{#1} (Restated).}\quad\itshape
}{%
  \par\medskip\normalfont
}

\numberwithin{equation}{section}
\newcommand{\spa}{\hspace{1pt}}
\newcommand{\vep}{\varepsilon}

\DeclarePairedDelimiter{\abs}{\lvert}{\rvert}

\newcommand{\norm}[1]{\left\lVert#1\right\rVert}
\newcommand{\setcond}{\hspace{2pt} \middle| \hspace{2pt}}

\DeclareFontFamily{U}{stix2bb}{\skewchar\font127 }
\DeclareFontShape{U}{stix2bb}{m}{n} {<-> stix2-mathbb}{}
\DeclareMathAlphabet{\mathbb}{U}{stix2bb}{m}{n}

\newcommand\lbar[1]{%
  \underaccent{\bar}{#1}}

\newcommand{\colA}[1]{\makebox[4.5em][c]{$#1$}}
\newcommand{\colB}[1]{\makebox[4.5em][c]{$#1$}}

\titleformat{\section}
  {\normalfont\bfseries\Large}     % font for section titles
  {\thesection}         % section number
  {1em}                 % space between number and title
  {\bfseries}        % makes the title bold
  
\titleformat{\subsection}
  {\normalfont\bfseries\normalsize}
  {\thesubsection}                
  {1em}                             
  {\bfseries}

\begin{document}

%%%title decoration%%%%%%%%%%%%%%%%

\newcommand\titlelowercase[1]{\texorpdfstring{\lowercase{#1}}{#1}}

%defines the size of author name
\font\mathptmx=cmr12 at 12pt

%%%%%%%%%%%%%%%%%%%%%%%%%%%%%%%

\title[\fontsize{13}{12}\mathptmx {\it{O\titlelowercase{n the intersection of} C\titlelowercase{antor sets and products of random matrices}}}]{\LARGE O\titlelowercase{n the intersection of} C\titlelowercase{antor sets} \protect{\\[9pt]} \titlelowercase{and products of random matrices}}

%available font sizes for title
%\tiny\scriptsize\footnotesize\small\normalsize\footnotesize\large\Large\LARGE\huge\Huge

\author[\fontsize{13}{12}\mathptmx {\it{N\titlelowercase{ima} A\titlelowercase{libabaei}}}]{\fontsize{13}{12}\mathptmx N\titlelowercase{ima} A\titlelowercase{libabaei}}



\maketitle

\begin{abstract}
Kenyon and Peres (1991) showed that the Hausdorff dimension of intersections of randomly translated Cantor sets can be expressed in terms of the top Lyapunov exponent of a product of random matrices, and this exponent can be written as an integral with respect to stationary measures on the projective line. Although explicit computations are available when stationary measures are discrete, the continuous case has remained challenging.

In this paper we introduce new combinatorial and analytic tools that allow us to compute the Lyapunov exponent, and hence the Hausdorff dimension, in a broad class of examples where stationary measures are continuous. As an application, we complete the dimension computation in the setting where a single digit is forbidden; for example, we determine the Hausdorff dimension of the intersection of the middle-seventh Cantor set with a random translate of itself.
\end{abstract}

\section{Introduction} \label{section: introduction}

\subsection{Background}

This paper studies the Hausdorff dimension of intersections of two Cantor-type sets under a random translation. While the constructions are elementary, the resulting dimension admits no simple closed form; instead, by Kenyon–Peres \cite{Kenyon--Peres: translated Cantor sets}, it is expressed in terms of the Lyapunov exponent of products of random nonnegative $2 \times 2$ matrices. The exponent can be written as an integral with respect to stationary measures, but it is typically difficult to compute explicitly. In this work, we introduce new techniques which allow us to explicitly calculate the Lyapunov exponent, and hence the dimension, in several interesting cases. As a consequence, the associated numerical scheme exhibits only polylogarithmic dependence on the target precision.

In what follows, we first review the results that motivated the present study, and our main contributions will be presented in the next subsection. Our starting point is an intriguing result due to Hawkes. Let $C$ be the classical, middle-third Cantor set. Hawkes proved the following result concerning the Hausdorff dimension of the intersection of translated Cantor sets.

\begin{theorem}[{\cite[Corollary of Theorem 1]{Hawkes}}] \label{theorem: Hawkes}
We have
\begin{equation*}
\mathrm{dim}_{\mathrm{H}} \left( C \cap (C + t) \right) = \frac{1}{3} \frac{\log{2}}{\log{3}} \, \, \, \, \left( \text{ Lebesgue a.e. $t \in [0, 1]$} \right).
\end{equation*}
\end{theorem}
For a while, this theorem stood as an isolated observation, until Kenyon and Peres revisited the question and brought it to a much deeper and more general level. Let $b \in \mathbb{Z}_{> 1}$. For $D_j \subset \{0, 1, \ldots, b-1\}$ ($j = 1, 2$), define
\begin{equation*}
K_j = \left\{ \sum_{n = 1}^\infty \frac{d_n}{b^n} \setcond d_n \in D_j \right\} \subset \mathbb{T} = \mathbb{R}/\mathbb{Z}.
\end{equation*}
Then, we are interested in $\mathrm{dim}_{\mathrm{H}} \left( (K_1 + t) \cap K_2 \right)$. Hawkes' method was not universal; essentially, it utilized the fact that the set $D_2 - D_1$ is contained in some arithmetic progression of length $b$. This is expanded as ``difference-set method'' within \cite[Section 3]{Kenyon--Peres: translated Cantor sets}.

Kenyon and Peres proved that, in general, the dimension in question can be expressed using the Lyapunov exponent of a product of some random matrices.
\begin{theorem}[{\cite[Theorem 1.2]{Kenyon--Peres: translated Cantor sets}}] \label{theorem: KP main result}
Let $\mu$ be a probability measure on $\mathbb{T} = \mathbb{R}/\mathbb{Z}$, invariant and ergodic for $T_b x = bx \mathrm{ \hspace{5pt} mod \hspace{3pt}} 1$. Define the $2$ by $2$ matrices $A_0, \ldots, A_{b-1}$ by
\begin{equation*}
A_u(i, j) = \# \big( \left( D_1 +i+u \right) \cap \left(D_2 + jb \right) \big) \quad \left( \text{ $0 \leq i, j \leq 1$ } \right),
\end{equation*}
where $\#$ is the number of elements. Then,
\begin{equation*}
\mathrm{dim}_{\mathrm{H}} \left( (K_1 + t) \cap K_2 \right)
= \frac{\lambda}{\log{b}}  \quad \left( \text{$\mu$-a.e. $t$} \right).
\end{equation*}
Here $\left( \text{for any choice of norm} \right)$
\begin{equation} \label{eq: lambda}
\lambda = \lim_{n \to \infty} \frac{1}{n} \log \left\| A_{t_n^{}} \cdots A_{t_1^{}} \right\| \text{\, $\bigg($ $\mu$-a.e. $t = \sum_{n = 1}^\infty \frac{t_n}{b^n}$ $\bigg)$}.
\end{equation}
Except when $A_{t_n^{}} \cdots A_{t_1^{}}$ is the zero matrix for $\mu$-a.e.\! $t$ for some $n$, in which case the dimension is $0$.
\end{theorem}

\begin{remark}
This $\lambda$ is called the (top) Lyapunov exponent of the product of random matrices $A_0, \ldots, A_{b-1}$. The RHS of \eqref{eq: lambda} may not always exist, but it exists and is constant for $\mu$-a.e. \@ $t$ (\cite{Furstenberg--Kesten}). Also, the original statement in \cite[Theorem 1.2]{Kenyon--Peres: translated Cantor sets} is formulated solely for the Lebesgue measure. However, a more general statement in \cite[Corollary 5.8]{Kenyon--Peres: translated Cantor sets} guarantees that this is true for any invariant and ergodic probability measure.
\end{remark}

Prior to Kenyon and Peres' work, Furstenberg and Kifer \cite{Furstenberg--Kifer} proved that the top Lyapunov exponent $\lambda$ can be expressed as an integral under the stationary measure on the projective space, as follows. Let $\norm{\cdot}_1$ denote the $\ell^1$-norm, the sum of the absolute values of entries.
The matrices $A_i$ are non-negative, so they induce maps $v \mapsto \frac{v \spa A_i}{ \norm{v \, A_i}_1}$ for a row vector $v$ in the simplex
\[ \Delta = \left\{ v = (v_1, v_2) \in \mathbb{R}^2 \setcond v_1, v_2 \geq 0,  v_1 + v_2 = 1 \right\}. \]
We identify $\Delta$ with $[-1, 1]$ with the map $[-1, 1] \ni x \mapsto \left( \frac{1+x}{2}, \frac{1-x}{2} \right) \in \Delta$. Then, the induced map above is identified with the following M\"obius transformation $f_i$ on $[-1,1]$.
\begin{equation} \label{eq: induced map by Ai}
f_i(x) =
\frac{
(p-q-r+s)x + (p-q+r-s) }{
(p+q-r-s)x + (p+q+r+s)}
\text{\quad for \quad}
A_i =
\begin{pmatrix}
p & q \\
r & s
\end{pmatrix}.
\end{equation}
This $\{f_i\}_{i=0}^{b-1}$ is said to be the \textbf{iterated function system (IFS)} associated with $(D_1, D_2)$.

A probability measure $\nu$ on $[-1,1]$ is said to be \textbf{$\boldsymbol{\mu}$-stationary} if for every continuous function $g$,
\[ \int_{[-1, 1]} g \, d\nu = \int_{[-1, 1]} \left( \int_{\mathbb{T}} g \circ f_{t_1^{}} (x) \, d\mu(t) \right) \, d\nu(x). \]
Then, \cite[Theorem 2.2]{Furstenberg--Kifer} states the following.
\begin{equation} \label{eq: lambda is the sup double int}
\lambda = \sup \left\{ \,
\int_{[-1, 1]} 
\left(
\int_{\mathbb{T}} \log \norm{
\begin{pmatrix}
\frac{1+x}{2}, & \hspace{-4pt} \frac{1-x}{2}
\end{pmatrix}
\, A_{t_1^{}} }_1 \, d\mu(t) \right)  \, d\nu(x)
\setcond \text{$\nu$ is $\mu$-stationary} \right\}.
\end{equation}
The uniqueness of $\mu$-stationary measure is guaranteed if
\begin{enumerate}
\item the semigroup generated by $\{A_i\}_i$ does not preserve any finite set of proper subspaces of $\mathbb{R}^2$, and one of $A_i$ (with $\mu(\{t | t_1 = i\}) > 0$) is primitive\footnote{A matrix $M$ is said to be primitive if there is a natural number $n$ such that every entry of $M^n$ is positive.}, or
\item one of $A_i$ (with $\mu(\{t | t_1 = i\}) > 0$) has rank one. (See \cite[Theorem 5.4]{Guivarch--Raugi}.)
\end{enumerate}

Using the equation \eqref{eq: lambda is the sup double int}, Kenyon and Peres showed that explicit calculation is possible in some cases. 

\begin{example}[{\cite[Example 1.3]{Kenyon--Peres: translated Cantor sets}}] \label{example: KP example}
Let $b = 4$ and $D_1 = D_2 = \{0,1,2\}$. We then have
\[
f_0(x) = 1, \quad f_1(x) = \frac{x+1}{2}, \quad f_2(x) = \frac{x-1}{2}, \quad f_3(x) = -1.
\]
Let $\mu$ be the Lebesgue measure, and $\nu$ be any $\mu$-stationary measure. The matrices $A_0$ and $A_3$ have rank one, so $\nu$ is unique and discrete (i.e. it is a countable sum of Dirac measures). The Lyapunov exponent can be calculated as
\[
\lambda = \frac{1}{6} \log{\frac{2}{3}} + \sum_{k=0}^\infty \frac{1}{4^{k+1}} \log{ \frac{(3 \cdot 2^k)!}{(2^{k+1})!} } = 0.797435 \cdots.
\]
This gives the dimension to be $ 0.575228 \cdots$.
\end{example}

The computation above exploits the discreteness of the stationary measure \cite[Remark, pp.~618]{Kenyon--Peres: translated Cantor sets}. Let us define the following.

\begin{definition}
The pair $(D_1, D_2)$ is \textbf{degenerate} if and only if there is a $\mu$-valid index $i$ such that $A_i$ has rank one.
\end{definition}
 When this happens, the image of the corresponding $f_i$ is a singleton, making the $\mu$-stationary measure unique and discrete. We can then calculate the dimension by tracking the trajectory of the said singleton under all $f_i$. The calculations by Kenyon and Peres (such as Example \ref{example: KP example}) fall under this category. See the recent paper by Fan and Vebitskiy \cite{Fan--Vebitskiy} for a general formula under this assumption.
 
In our setting, the formula in \cite[Theorem 1.1]{Fan--Vebitskiy} writes as follows. Consider the same setting as in Theorem \ref{theorem: KP main result}, and let $\{f_i\}_{i \in I}$ be the associated IFS where $I = \{0, 1, \ldots, b-1\}$. Suppose that there is a $\mu$-valid index $j$ and $\alpha \in [-1, 1]$ with $f_j([-1, 1]) = \{ \alpha \}$. Then, letting $p_i = \mu( \{ t | t_1 = i \} )$ for each $i$,
\begin{equation*}
\lambda = p_j \sum_{k \in I} p_k \sum_{ \lbar{i} \in ( I \backslash \{j\} )^* } p_{\lbar{i}} \log{ \norm{
\begin{pmatrix}
\frac{1+ f_{ \lbar{i} } (\alpha) }{2}, & \hspace{-4pt} \frac{1- f_{ \lbar{i} } (\alpha) }{2}
\end{pmatrix} A_k }_1 }. 
\end{equation*}
Here, the notations are $( I \backslash \{j\} )^* = \{ \varnothing \} \cup \left( \bigcup_{n=1}^\infty ( I \backslash \{j\} )^n \right)$, $f_{ \lbar{i} } = f_{i_n} \circ \cdots \circ f_{i_1}$, $p_{ \lbar{i} } = p_{i_1} \cdots p_{i_n}$ for $\lbar{i} = (i_1, \ldots, i_n)$, and $f_{\varnothing}$ is the identity map. This is immediate, upon seeing that for a continuous function $g$ on $[-1, 1]$,
\[ \int_{[-1, 1]} g \, d\nu = p_j \, g \circ f_j( \alpha) + \sum_{i \in ( I \backslash \{j\} )} p_i \int_{[-1, 1]} g \circ f_i \spa d\nu. \]

Consequently, the computation for degenerate cases is already well-understood. Beyond this setting, where the stationary measure is continuous, one confronts fundamentally different obstacles, which we aim to tackle in this study.

One important line of research was pioneered by Ruelle \cite{Ruelle76} and Pollicott \cite{Pollicott10}, and later improved by Jurga-Morris \cite{Jurga-Morris}. They developed a spectral–analytic approach in which the Lyapunov exponent is expressed in terms of the Fredholm determinant of an associated transfer operator, thereby reducing its computation to finite combinatorial data. The method applies to a wide class of random matrix products; however, the computational cost grows subexponentially. To achieve a truncation error $\leq \vep$, the cost is
\[
O \! \left( \exp{ \left(  \sqrt{ \log{ \left( 1/ \vep \right) } } \right)  } \right).
\]

The methods we develop in this paper are not applicable to all random matrix products, but whenever the required hypotheses hold, we obtain the same accuracy with polylogarithmic complexity (see Remark \ref{remark: computational cost}):
\[
O \! \left( \Big( \log{ \left( 1/\vep \right) } \Big)^3  \right).
\]
For instance, we will see that our new method completes the calculation of the Hausdorff dimension of $ \left( K_1 + t \right) \cap K_2 $ in all cases with $b \geq 7$ in which exactly one digit is forbidden from each $D_j$.

\, \\[-25pt]

\subsection{Main results with examples} \label{subsection: Main results with examples}
We provide an explicit algorithm to compute the Hausdorff dimension of the intersection of translated Cantor sets in several cases. In fact, our method yields two distinct scenarios in which the Hausdorff dimension can be computed explicitly.

First is when $\{f_i\}_i$ from \eqref{eq: induced map by Ai} admits a common structure in a sense.

\begin{definition}
For a probability measure $\mu$ on $\mathbb{T}$, an index $i \in \{ 0, 1, \ldots, b-1\}$ is \textbf{$ \boldsymbol{\mu}$-valid} if
\[
\mu( t_1 = i ) := \mu \left( \left\{ t \in \mathbb{T} \setcond t_1 = i \right\} \right) > 0.
\]
A continuous function $\phi: [-1, 1] \to [-1, 1]$ is said to be \textbf{affine co-invariant} with respect to $\mu$ and $\{f_i\}_i$ if for each $\mu$-valid $0 \leq i \leq b-1$, there are $r_i, s_i \in \mathbb{R}$ such that
\[ \phi \circ f_i = r_i \phi + s_i. \]
\end{definition}
For example, if $\{f_i\}_i$ are all affine maps, the function $\phi(x) = x$ is affine co-invariant. If such a function exists and can be used to expand the integrand in \eqref{eq: lambda is the sup double int}, the dimension admits a closed form expression (Proposition \ref{proposition: affine co-invariant method}). Here, we will see an example. 

Let $I = \{0, 1, \ldots, b-1\}$.  Let $\tilde{\mu}$ be the probability measure on $I$ defined by a probability vector $\{p_i\}_{i \in I}$, that is, $\tilde{\mu}(\{i\}) = p_i \geq 0$ and $\sum_i p_i = 1$. Consider the map $\pi: I^\mathbb{N} \to \mathbb{T}$ defined by
\[
\pi: I^\mathbb{N} \ni \left( d_1, d_2, \ldots \right) \mapsto \sum_{n=1}^\infty \frac{d_n}{b^n} \in \mathbb{T}.
\]
Let $\mu$ be the push-forward of the product measure $\tilde{\mu}^{\mathbb{N}}$ on $I^{\mathbb{N}}$ by $\pi$.
We abuse the notation and denote as $\mu = \left( p_0, p_1, \cdots, p_{b-1} \right)^{\otimes \mathbb{N}}$ and call this a \textbf{product measure}.

\begin{example} \label{ex: recurring method}
Let $b = 7, D_1 = \{0, 2, 5\},$ and $D_2 = \{0, 1, 2, 4, 6\}$. Let $\phi(x) = \frac{6x+1}{-4x+12}$. Then, $\phi$ is affine co-invariant under the following maps.
\[
f_1(x) = \frac{x+3}{-x+5}, \quad f_3(x) = \frac{x}{x-2}.
\]
Let $\mu = \left( 0, \frac{1}{2}, 0, \frac{1}{2}, 0, 0, 0 \right)^{\otimes \mathbb{N}}$. Note that $f_1$ and $f_3$ have distinct attracting fixed point, $1$ and $0$ respectively, so the semigroup generated by $\{A_i\}_i$ does not preserve any finite set of proper subspaces of $\mathbb{R}^2$. This implies the uniqueness of $\mu$-stationary measure $\nu$. Moreover, $\nu$ is continuous by \cite[Lemma 4.2]{Bougerol--Lacroix}, so the method by Kenyon and Peres does not work. Take the $\mu$-stationary measure $\nu$ and let $\xi_n = \int \phi^n d\nu$. Then, it satisfies the following recurrence formula.
\[
\xi_n = \left( 1 - \frac{ 1 + (-1)^n }{ 2^{n+1} } \right)^{-1} \sum_{k=0}^{n-1} \binom{n}{k} \left( \left( \frac{7}{16} \right)^{n-k} + (-1)^k \left( \frac{1}{8} \right)^{n-k} \right) \frac{ \xi_k }{ 2^{k+1} }.
\]
(It has an explicit expression as well. See Lemma \ref{lemma: calculation of integral of co-affine function} for details.) Then, we have
\begin{align*}
\lambda
&= \frac{1}{2} \left( \log{13} + \log{31} - \log{72} \right) + \sum_{n=1}^{\infty} \frac{1}{2n} \left\{ (-1)^{n-1} \left( \left( \frac{8}{31} \right)^n - 2 \cdot \left( \frac{2}{3} \right)^n \right) - \left( \frac{4}{13} \right)^n \right\} \xi_n \\
&= 0.693147\cdots.
\end{align*}
Therefore,
\begin{equation*}
\mathrm{dim}_{\mathrm{H}} \left( (K_1 + t) \cap K_2 \right)
= \frac{ \lambda }{ \log{7} } = 0.356207\cdots \quad \left( \text{$\mu$-a.e. $t$} \right).
\end{equation*}
\end{example}

Next, in a more analytic vein, we introduce the notion of Kernel-expansion to characterize pairs of Cantor-type sets for which the integrand in \eqref{eq: lambda is the sup double int} can be expanded in terms of kernel functions, which integrate to zero. In practice, this decomposition is carried out using Neumann series, which provides an explicit representation of the integral and thus allows us to compute the Hausdorff dimension.

Let $\mu$ be a probability measure on $\mathbb{T} = \mathbb{R}/\mathbb{Z}$, invariant and ergodic for $T_b x = bx \mathrm{ \hspace{5pt} mod \hspace{3pt}} 1$. For a continuous function $\phi: [-1, 1] \to (-1, 1)$, let
\begin{equation*}
\psi_n(\phi) = \phi^n \, - \, \int_{\mathbb{T}} \left(\phi \circ f_{t_1^{}} \right)^n d\mu(t). \quad \text{( $n \in \mathbb{N}$ )}
\end{equation*}
Note that due to the $\mu$-stationarity (or, self-similarity) of $\nu$, we have
\[ \int \psi_n(\phi) \, d\nu = 0. \]
Let $\mathbb{N}_0 = \{0, 1, 2, \ldots\}$, and $\ell^\infty(\mathbb{N}_0)$ be the space of bounded sequences in $\mathbb{R}$ indexed with $\mathbb{N}_0$.
\begin{definition}
The pair $(D_1, D_2)$ (or the matrices $\{A_i\}_i$) is said to be \textbf{Kernel-expandable with respect to $\boldsymbol{\mu}$} if there is a continuous function $\phi: [-1, 1] \to (-1, 1)$ such that for every $\mu$-valid $0 \leq i \leq b-1$, there is $(a_n)_{n \in \mathbb{N}_0} \in \ell^\infty(\mathbb{N}_0)$ satisfying
\begin{equation} \label{eq: definition of kernel expandable}
\log \norm{
\begin{pmatrix}
\frac{1+x}{2}, & \hspace{-4pt} \frac{1-x}{2}
\end{pmatrix}
\, A_i }_1
= a_0 + \sum_{n=1}^\infty a_n \psi_n(\phi(x)).
\end{equation}
When $\mu$ is the Lebesgue measure we will simply say that $(D_1, D_2)$ is Kernel-expandable.
\end{definition}

Note that by $\| \psi_n(\phi) \|_\infty \leq 2 {\| \phi \|_\infty}^n$ and $\| \phi \|_\infty < 1$, the uniform convergence of the right hand side in the above equation is always ensured.

Informally, Kernel-expansion decomposes the integrand into a constant plus a uniformly convergent series of functions each having zero $\nu$-average. We remark that if $\{A_i\}_i$ is Kernel-expandable, the integral in equation \eqref{eq: lambda is the sup double int} takes the same value for every $\mu$-stationary measure. In other words, if there are $\mu$-stationary measures that yield different values of this integral, then $\{A_i\}_i$ is not Kernel-expandable.

Theorem \ref{theorem: NAC assures KE} states that if we can find a $\phi$ that behaves well, we can calculate the top Lyapunov exponent of the product of random matrices, and thus the dimension. Here, we will state a rough version for conciseness.

\begin{theorem*} [Rough version of Theorem \ref{theorem: NAC assures KE}]
Let $\mu$ be a probability measure on $\mathbb{T} = \mathbb{R}/\mathbb{Z}$, invariant and ergodic for $T_b x = bx \mathrm{ \hspace{5pt} mod \hspace{3pt}} 1$. Define an infinite matrix $T = (b_{k, n})_{k, n \in \mathbb{N}_0}$ by $b_{k,0} = 0$ for all $k$, and for $n > 0$, as the coefficients appearing in the following expansion (where we assume such expansion exists).
\[
\sum_{i = 0}^{b-1} \, \mu(t_1 = i) \,  ( \phi \circ f_i )^n = \sum_{k=0}^\infty b_{k,n} \phi^k.
\]
If we can find a ``nice'' $\phi$ such that
\begin{enumerate}
\item $T: \ell^\infty(\mathbb{N}_0) \to \ell^\infty(\mathbb{N}_0)$ is well-defined and $\norm{T} < 1$.
\item For each $\mu$-valid index $i$, we can find $a^{(i)} = \left( a_n^{(i)} \right)_{n=0}^\infty \in \ell^\infty(\mathbb{N}_0)$ satisfying the following equality with uniformly convergent series.
\begin{equation*}
\log \norm{
\begin{pmatrix}
\frac{1+x}{2}, & \hspace{-4pt} \frac{1-x}{2}
\end{pmatrix}
\, A_i }_1
=
\sum_{n=0}^\infty a_n^{(i)} \phi(x)^n  \quad \left(\, x \in [-1, 1] \, \right).
\end{equation*}
\end{enumerate}
Then, the pair $(D_1, D_2)$ is Kernel-expandable with respect to $\mu$, and
\begin{equation}
\mathrm{dim}_{\mathrm{H}} \left( (K_1 + t) \cap K_2 \right)
=
\frac{1}{\log{b}} \, \left( \sum_{n=0}^\infty T^n \, a \right)_0 \quad \left( \text{$\mu$-a.e. $t$} \right),
\end{equation}
where $a = \sum_{i=0}^{b-1}  \, \mu(t_1 = i) \,   a^{(i)}$. $\left( \text{Here, $(\cdot)_0$ is the $0$-th entry.} \right)$
\end{theorem*}

Theorem \ref{theorem: NAC assures KE} is unique in that it tackles the difficult case of when the $\mu$-stationary measure is continuous. The condition imposed on $\phi$ is called the \textbf{Neumann Admissibility Condition} (or NAC), and will be defined in equation \eqref{condition: NAC} in \S\S\ref{subsection: introducing the method}. A natural question is whether such a ``good'' $\phi$ can actually be found. In our setting, the answer is affirmative in several interesting cases. The next theorem ensures Kernel-expandability by finding a $\phi$ satisfying NAC.

\begin{theorem} \label{theorem: toothless twins}
Let $b \geq 7$ be a natural number and $\mu$ Lebesgue measure. Consider $\tau, u \in \{0, 1, \ldots, b-1\}$ and let
\begin{align*}
& D_1 = \{0, 1, \ldots, b-1\} \backslash \{ \tau \}, \\
& D_2 = \{0, 1, \ldots, b-1\} \backslash \{u\}.
\end{align*}
Then, the pair $(D_1, D_2)$ is either degenerate or Kernel-expandable. 
More precisely, it is Kernel-expandable if both $\tau \notin \{0, b-1\}$ and $\tau + u = b-1$. Otherwise, it is degenerate.
\end{theorem}

This theorem shows that when $b \geq 7$ and we forbid only a single digit for each $D_j$, we have completed the calculation of $\mathrm{dim}_{\mathrm{H}} \left( (K_1 + t) \cap K_2 \right)$ for Lebesgue a.e.\! $t$. That is, either the already-known method by Kenyon and Peres can be used, or our new technique is viable.

Let us see an example where Kernel-expansion yields the dimension for a non-degenerate case.

\begin{example} \label{example: middle seventh}
Let $b = 7$, and $D_1 = D_2 = \{0, 1, 2, 4, 5, 6\}$. So, we are considering the middle-seventh Cantor sets. Let $\mu$ be the Lebesgue measure. Then, the matrices $A_0, \ldots, A_6$ are, by order of appearance,
\[
\begin{pmatrix}
 6 & 0\\
 4 & 1
\end{pmatrix}
, \hspace{2pt}
\begin{pmatrix}
 4 & 1\\
 3 & 2
\end{pmatrix}
, \hspace{2pt}
\begin{pmatrix}
3 & 2\\
2 & 3
\end{pmatrix}
, \hspace{2pt}
\begin{pmatrix}
2 & 3\\
3 & 2
\end{pmatrix}
, \hspace{2pt}
\begin{pmatrix}
3 & 2\\
2 & 3
 \end{pmatrix}
, \hspace{2pt}
\begin{pmatrix}
2 & 3\\
1 & 4
\end{pmatrix}
, \hspace{2pt}
\begin{pmatrix}
1 & 4\\
0 & 6
\end{pmatrix}.
\]
The associated IFS $\{f_i\}_i$ are
\[
f_0(x) = \frac{3x+9}{x+11}, \hspace{3pt}
f_1(x) = \frac{x+2}{5}, \hspace{3pt}
f_2(x) = \frac{x}{5}, \hspace{3pt}
f_3(x) = -\frac{x}{5},
\]
\[
f_4(x) = \frac{x}{5}, \hspace{3pt}
f_5(x) = \frac{x-2}{5}, \hspace{3pt}
f_6(x) = \frac{3x-9}{-x+11}.
\]
Since $f_0$ and $f_6$ have distinct attracting fixed points, the semigroup generated by $\{A_i\}_i$ does not preserve any finite set of proper subspaces of $\mathbb{R}^2$, implying the uniqueness of $\mu$-stationary measure $\nu$. Furthermore, $\nu$ is continuous by \cite[Lemma 4.2]{Bougerol--Lacroix}. However, Kernel expansion is possible with $\phi(x) = \frac{7x+4}{4x+16}$. The operator norm of $T$ in the above ``Theorem'' is less than $\frac{649}{723}$ for this $\phi$. Define $a = (a_n)_{n=0}^\infty \in \ell^\infty(\mathbb{N}_0)$ by
\begin{align*}
\hspace{60pt} &a_0 = 2 \log{ \frac{11}{2} } + 5\log{5}+ \log{ \frac{73}{77} } + \log{ \frac{81}{77} },& \\
&a_k = \frac{1}{k} \left\{ \left( \frac{28}{73} \right)^k - 2 \left( \frac{4}{7} \right)^k + \left( \frac{60}{81} \right)^k \right\} \quad \text{for $k > 0$}. &
\end{align*}
We can then calculate $\lambda$ as \\[-15pt]
\begin{align*}
 \lambda = \sum_{n=0}^\infty \left( T^n a \right)_0 = 1.6363797884 \cdots.
\end{align*}
Therefore, \\[-20pt]
\begin{equation*}
\mathrm{dim}_{\mathrm{H}} \left( (K_1 + t) \cap K_2 \right) = \frac{\lambda}{\log{7}} = 0.8409328607 \cdots.
\end{equation*}
The function $\phi$ used above is convenient but not optimal for numerical convergence; the natural numbers $M, N$ appearing in Proposition \ref{proposition: evaluation of Neumann series} are $112$ and $697$ respectively, for the tolerated error of $\vep = 10^{-4}$. For this example, $\phi_0(x) = \frac{0.3769x -0.2768}{-0.1973x+1}$ gives faster convergence, when we have $(M, N) = (87, 182)$ for the same $\vep$. For $\vep = 10^{-50}$, we have $(M, N) = (863, 1251)$ under $\phi_0$. (Note that we have $M, N \sim \log(1/\vep)$ by Proposition \ref{proposition: evaluation of Neumann series}.) 
\end{example}

The kernel-expandability is not limited to the case where we forbid a single digit.

\begin{example}
Let $b = 9$, $\mu$ Lebesgue measure, and
\[ D_1 = \{ 0, 1, 4, 5, 7, 8 \} = \{ 0, 1, \ldots, 8 \} \backslash \{ 2, 3, 6 \}, \]
\[ D_2 = \{ 0, 2, 3, 5, 6, 8\} = \{ 0, 1, \ldots, 8 \} \backslash \{ 1, 4, 7 \}. \]
The matrices $A_0, \ldots, A_8$ are, by order of appearance,
\begin{samepage}
\[
\begin{pmatrix}
3 & 0\\
4 & 1
\end{pmatrix}
,\hspace{6pt}
\begin{pmatrix}
4 & 1\\
3 & 1
\end{pmatrix}
,\hspace{6pt}
\begin{pmatrix}
3 & 1\\
2 & 1
\end{pmatrix}
,\hspace{6pt}
\begin{pmatrix}
2 & 1\\
2 & 3
\end{pmatrix}
,\hspace{6pt}
\begin{pmatrix}
2 & 3\\
2 & 2
\end{pmatrix}
,
\]
\[
\begin{pmatrix}
2 & 2\\
1 & 2
\end{pmatrix}
,\hspace{6pt}
\begin{pmatrix}
1 & 2\\
1 & 4
\end{pmatrix}
,\hspace{6pt}
\begin{pmatrix}
1 & 4\\
1 & 3
\end{pmatrix}
,\hspace{6pt}
\begin{pmatrix}
1 & 3\\
0 & 3
\end{pmatrix}.
\]
\end{samepage}

The associated IFS $\{f_i\}_i$ are
\[
f_0(x)=\frac{3}{-x+4},\hspace{10pt}
f_1(x)=\frac{x+5}{x+9},\hspace{10pt}
f_2(x)=\frac{x+3}{x+7},
\]
\[
f_3(x)=\frac{x}{-x+4},\hspace{10pt}
f_4(x)=\frac{-x-1}{x+9},\hspace{10pt}
f_5(x)=\frac{x-1}{x+7},
\]
\[
f_6(x)=\frac{x-2}{-x+4},\hspace{10pt}
f_7(x)=\frac{-x-5}{x+9},\hspace{10pt}
f_8(x)=\frac{x-5}{x+7}.
\]
So, in this case, none of the $f_i$ is linear, and we do not see the type of symmetry in matrices $A_i$ that is present in the previous Example \ref{example: middle seventh}. A similar discussion to that in Example \ref{example: middle seventh} shows that the $\mu$-stationary measure is unique and continuous; yet $(D_1, D_2)$ is Kernel-expandable with $\phi(x) = \frac{-0.3914x -0.055}{-0.0639x +1}$, and we have
\[
\lambda = 1.3770228916\cdots.
\]
Thus,
\begin{equation*}
\mathrm{dim}_{\mathrm{H}} \left( (K_1 + t) \cap K_2 \right) = \frac{\lambda}{\log{9}} =0.6267101259\cdots.
\end{equation*}
\end{example}

From numerical experiments, we conjecture that if the number of forbidden digits is relatively small, then the digits pair is either degenerate or Kernel-expandable.
\begin{conjecture}
There is a number $0 < r < 1$ close to $1$ such that for a large enough $b$, if
\[ \frac{ \# D_1 }{ b } > r, \text{\, and \,} \frac{ \# D_2 }{ b } > r, \]
then the pair $(D_1, D_2)$ is either degenerate or Kernel-expandable. 
\end{conjecture}

\section{Recurring method}

Here we state the closed form expression for the Hausdorff dimension of $(K_1 + t) \cap K_2$ when we have a nice affine co-invariant function.

Take $b \in \mathbb{Z}_{\geq 2}$ and let $I = \{0, 1, \ldots, b-1\}$. Consider $D_1, D_2 \subset I$ and let $\{A_i\}_{i \in I}$ be the matrices in Theorem \ref{theorem: KP main result}. Recall the following M\"obius transformation $f_i$ on $[-1,1]$ induced by $A_i$ in \eqref{eq: induced map by Ai}.
\begin{equation*}
f_i(x) =
\frac{
(p-q-r+s)x + (p-q+r-s) }{
(p+q-r-s)x + (p+q+r+s)}
\text{\quad for \quad}
A_i =
\begin{pmatrix}
p & q \\
r & s
\end{pmatrix}.
\end{equation*}
This $\{f_i\}_{i \in I}$ was called the IFS associated with $(D_1, D_2)$. (These $f_i$ may map all elements to a single point, which happens when $A_i$ has rank $1$.)

\begin{lemma} \label{lemma: calculation of integral of co-affine function}
Suppose there is an affine co-invariant function $\phi: [-1,1] \to [-1,1]$ with respect to $\mu$ and $\{f_i\}_i$. For each $\mu$-valid $0 \leq i \leq b-1$, let $r_i, s_i \in \mathbb{R}$ be the real constants such that
\[ \phi \circ f_i = r_i \phi + s_i. \]
Let $p_i = \mu( \{ t | t_1 = i \} )$ for each $i$. For natural numbers $k$ and $n$, define a positive number $\Lambda_{k, n}$ by
\[ \Lambda_{k, n} = \frac{\sum_{i \in I} p_i \, s_i^n \, r_i^k}{1 - \sum_{i \in I} p_i \, r_i^k}. \]
Then, for any $\mu$-stationary measure $\nu$ and a natural number $n$,
\begin{equation} \label{eq: integral of phi}
\int_{[-1,1]} \, \phi^n \,d\nu = 
\frac{1}{1 - \sum_{i \in I} p_i r_i^n} \hspace{2pt}
\sum_{ \substack{ n_1 + \cdots + n_k = n \\
 k \in \mathbb{N}, n_\ell \geq 1}} \hspace{2pt}
\frac{n!}{n_1! \cdots n_k!} \, \Lambda_{0, n_1} \Lambda_{n_1, n_2} \cdots \Lambda_{n_{k-1}, n_k}.
\end{equation}
\end{lemma}

\begin{breakremark}
\noindent (1) Let us denote by $\xi_n$ the integral of $\phi^n$ with respect to a $\mu$-stationary measure $\nu$. It has the following recurrence relation:
\begin{equation} \label{eq: recurrence of xi}
\xi_n = \left( 1 - \sum_{i \in I} p_i r_i^n \right)^{-1} \hspace{3.5pt} \sum_{k = 0}^{n-1}
\begin{pmatrix}
n \\
k
\end{pmatrix}
\left(
\sum_{i \in I} p_i \spa r_i^k \spa s_i^{n-k}
\right) \xi_k \quad \left( \text{ for $n \geq 1$ } \right).
\end{equation}
For computations it is considerably faster to use this recurrence than the explicit combinatorial expansion in equation \eqref{eq: integral of phi}. \\
\noindent (2) The idea of using self-similarity to calculate the integral of polynomials appears in \cite{Strichartz}.
\end{breakremark}

\begin{proof}
Let $\nu$ be a $\mu$-stationary measure. Let
\[ \xi_n = \int_{[-1,1]} \, \phi^n \, d\nu. \]
We have $\xi_0 = 1$. By the $\mu$-stationarity (self-similarity) of $\nu$,
\begin{align*}
\xi_n
&= \sum_{i \in I} p_i \int_{[-1,1]} \, \left( \phi \circ f_i \right)^n \, d\nu \\
&= \sum_i p_i \int_{[-1,1]} \big( r_i \phi + s_i \big)^n \spa d\nu \\
&= \sum_i p_i \spa \sum_{k = 0}^n
\begin{pmatrix}
n \\
k
\end{pmatrix}
r_i^k \spa s_i^{n-k} \spa \xi_k,
\end{align*}
arriving at the recurrence relation in equation \eqref{eq: recurrence of xi}. The equation \eqref{eq: integral of phi} follows by induction by iterating the equation \eqref{eq: recurrence of xi}.
\end{proof}

\begin{proposition} \label{proposition: affine co-invariant method}
Suppose there is an affine co-invariant function $\phi: [-1,1] \to [-1,1]$ with respect to $\mu$ and $\{f_i\}_i$ such that for each $\mu$-valid $0 \leq i \leq b-1$, there is $\left( a_n^{(i)} \right)_n \in \mathbb{R}^{\mathbb{N}_0}$ satisfying the following equality with uniformly convergent series.
\begin{equation*}
\log \norm{
\begin{pmatrix}
\frac{1+x}{2}, & \hspace{-4pt} \frac{1-x}{2}
\end{pmatrix}
\, A_i }_1
=
\sum_{n=0}^\infty a_n^{(i)} \phi(x)^n  \quad \left(\, x \in [-1, 1] \, \right).
\end{equation*}
Then, the Hausdorff dimension of $(K_1 + t) \cap K_2$ is calculated under the notation of Lemma \ref{lemma: calculation of integral of co-affine function} as
\begin{flalign*}
& \mathrm{dim}_{\mathrm{H}} \left( (K_1 + t) \cap K_2 \right)
&
\end{flalign*}
\begin{flalign*}
\hspace{15pt} = \frac{1}{\log{b}} \,
\sum_{n = 0}^{\infty}
\left( \sum_{i \in I} p_i \, \frac{ a_n^{(i)} }{1 - \sum_{i \in I} p_i r_i^n} \hspace{2pt} \right)
\sum_{ \substack{ n_1 + \cdots + n_k = n \\
 k \in \mathbb{N}, n_\ell \geq 1}} \hspace{2pt}
\frac{n!}{n_1! \cdots n_k!} \, \Lambda_{0, n_1} \Lambda_{n_1, n_2} \cdots \Lambda_{n_{k-1}, n_k}
\quad \left( \text{$\mu$-a.e. $t$} \right).
\end{flalign*}
\end{proposition}

\begin{proof}
By Lemma \ref{lemma: calculation of integral of co-affine function}, the integral of $\phi^n$ with respect to any $\mu$-stationary measure is calculated as in equation \eqref{eq: integral of phi}. Since the expansion with $\left( a_n^{(i)} \right)_n$ is uniformly convergent, and by  equation \eqref{eq: lambda is the sup double int}, we obtain the desired expression.
\end{proof}

This proposition is used to carry out the calculation in Example \ref{ex: recurring method}

\section{Neumann series method}

\subsection{Introducing the method} \label{subsection: introducing the method}

In this section, we introduce the most novel idea in this paper which uses Neumann series. We take a M\"obius transformation $\phi$ and study the expansion of each $\phi \circ f_i$ in the $\phi$-basis. For application, this subsection is written in terms of general finite number of $2 \times 2$ matrices with non-negative entries. Also, the result in this section (Proposition \ref{prop: KE core}) can be viewed as a preliminary form of Theorem \ref{theorem: NAC assures KE}, which we will prove in the next subsection. 

We will employ the following notation, where we always assume $\gamma x + \delta \ne 0$ for $x \in [-1, 1]$.
\[
\begin{pmatrix}
\alpha & \beta \\
\gamma & \delta
\end{pmatrix}(x)
=
\frac{\alpha x + \beta}{\gamma x + \delta}.
\]

Let $b \geq 2$ be an integer, and $\mu$ a probability measure on $\mathbb{T} = \mathbb{R}/\mathbb{Z}$, invariant and ergodic for $T_b x = bx \mathrm{ \hspace{5pt} mod \hspace{3pt}} 1$. Suppose we have $b$ number of non-zero $2 \times 2$ matrices with non-negative entries: $A_0, \ldots, A_{b-1} \in M_2(\mathbb{R}_{\geq 0})$. Let $\{f_i\}_{i=0}^{b-1}$ be the associated M\"obius IFS; i.e. $f_i: [-1, 1] \to [-1, 1]$ is defined with
\begin{equation*}
f_i(x) =
\begin{pmatrix}
p-q-r+s & p-q+r-s \\
p+q-r-s & p+q+r+s
\end{pmatrix}(x)
\text{\quad for \quad}
A_i =
\begin{pmatrix}
p & q \\
r & s
\end{pmatrix}.
\end{equation*}

The following proposition is the core of Neumann series method.
\begin{proposition} \label{prop: KE core}
Let $\phi: [-1, 1] \to (-1, 1)$ be a M\"obius transformation, and $\mu$ a probability measure on $\mathbb{T} = \mathbb{R}/\mathbb{Z}$, invariant and ergodic for $T_b x = bx \mathrm{ \hspace{5pt} mod \hspace{3pt}} 1$. Suppose we have the following expansion for every natural number $n$:
\[
\sum_{i = 0}^{b-1} \mu( t_1 = i ) \, ( \phi \circ f_i )^n = \sum_{k=0}^\infty b_{k,n} \phi^k.
\]
Set $b_{k, 0} = 0$ for all $k \in \mathbb{N}_0$, and define an infinite matrix $T = (b_{k, n})_{k, n \in \mathbb{N}_0}$. We assume that
\begin{enumerate}
\item $T: \ell^\infty(\mathbb{N}_0) \to \ell^\infty(\mathbb{N}_0)$ is well-defined and $\norm{T} < 1$.
\item For each $\mu$-valid index $i$, we can find $a^{(i)} = \left( a_n^{(i)} \right)_{n=0}^\infty \in \ell^\infty(\mathbb{N}_0)$ satisfying the following equality with uniformly convergent series.
\begin{equation*}
\log \norm{
\begin{pmatrix}
\frac{1+x}{2}, & \hspace{-4pt} \frac{1-x}{2}
\end{pmatrix}
\, A_i }_1
=
\sum_{n=0}^\infty a_n^{(i)} \phi(x)^n  \quad \left(\, x \in [-1, 1] \, \right).
\end{equation*}
\end{enumerate}
Then, $\{A_i\}_{i=0}^{b-1}$ is Kernel-expandable, and its top Lyapunov exponent $\lambda$ can be calculated as
\begin{equation*}
\lambda
=
\left( \sum_{n=0}^\infty T^n \, a \right)_0,
\end{equation*}
where $a = \sum_{i=0}^{b-1} \, \mu( t_1 = i ) \, a^{(i)}$.
\end{proposition}

\begin{proof}

For $n \in \mathbb{N}$, define
\begin{equation*}
\psi_n = \phi^n \, - \, \sum_{i = 0}^{b-1} \mu( t_1 = i ) \, ( \phi \circ f_i )^n.
\end{equation*}
These $\psi_n$ are ``kernel'' functions --- for any $\mu$-stationary measure, $\nu$, we have
\[
\int_{[-1, 1]} \psi_n d\nu = 0.
\]

Consider the following function spaces.
\begin{equation*}
\mathscr{N} := \left\{ h: [-1, 1] \to \mathbb{R}
\setcond
\text{There is $c = (c_n)_n \in \ell^\infty(\mathbb{N}_0)$ with } h = c_0 + \sum_{n=1}^\infty c_n \psi_n \right\},
\end{equation*}
\begin{equation*}
\mathscr{A} := \left\{ g: [-1, 1] \to \mathbb{R}
\setcond
\text{There is $ a =(a_k)_k \in \ell^\infty(\mathbb{N}_0)$ with } g = \sum_{k=0}^\infty a_k \phi^k \right\}.
\end{equation*}
By the assumption that $\mathrm{Im}(\phi) \subset (-1, 1)$, we have $\| \phi \|_\infty < 1$. Since $\| \psi_n(\phi) \|_\infty \leq 2 {\| \phi \|_\infty}^n$, the series inside the definitions above are uniformly convergent.

Let $I: \ell^\infty(\mathbb{N}_0) \to \ell^\infty(\mathbb{N}_0)$ be the identity map. By the definition of $\psi_n$ and $T$, the following composition of operators is an inclusion.
\[
\mathscr{N} \xrightarrow{\cong} \ell^\infty(\mathbb{N}_0) \xrightarrow{I - T} \ell^\infty(\mathbb{N}_0) \xrightarrow{\cong} \mathscr{A},
\]
where the isomorphism  $\mathscr{N} \cong \ell^\infty(\mathbb{N}_0)$ and $\mathscr{A} \cong \ell^\infty(\mathbb{N}_0)$ are by the correspondence of $c$ and $a$, respectively. Since $\norm{T} < 1$, this $I - T$ is invertible using Neumann series:
\begin{equation*}
\left( I - T \right)^{-1} = \sum_{n=0}^\infty T^n.
\end{equation*}

Now, let $h = c_0 + \sum_{n=1}^\infty c_n \psi_n \in \mathscr{N}$. We then have
\begin{equation*}
\int_{[-1, 1]} h \, d\nu = c_0.
\end{equation*}
This implies that, for any $g =  \sum_{n=0}^\infty u_k \phi^k \in \mathscr{N}$ with $u = \left( u_k \right)_k$,
\begin{align*}
\int_{[-1, 1]} g \, d\nu \,
&= \, \int_{[-1, 1]}
\begin{pmatrix}
1, \, \psi_1, \, \psi_2, \, \cdots
\end{pmatrix}
\sum_{n=0}^\infty T^n \, u \, d\nu \, \\
&= \left( \sum_{n=0}^\infty T^n \, u \right)_0.
\end{align*}
Here, we used the uniform convergence of the series due to $\norm{T} < 1$. We assumed that the integrand in the following expression (repeat of equation \eqref{eq: lambda is the sup double int}) has a uniformly convergent expansion in terms of $\phi$.
\begin{equation*}
\lambda = \sup \left\{ \,
\int_{[-1, 1]} \int_{\mathbb{T}} \log \norm{
\begin{pmatrix}
\frac{1+x}{2}, & \hspace{-4pt} \frac{1-x}{2}
\end{pmatrix}
\, A_{t_1^{}} }_1 \, d\mu(t) \, d\nu(x)
\setcond \text{$\nu$ is $\mu$-stationary} \right\}.
\end{equation*}
This implies
\begin{equation*}
\lambda
=
\left( \sum_{n=0}^\infty T^n \, a \right)_0,
\end{equation*}
letting $a = \sum_{i=0}^{b-1} \, \mu( t_1 = i ) \, a^{(i)}$.
\end{proof}

\subsection{Neumann Admissibility Condition and Kernel-expansion}

In this subsection, we consider an analytic condition for $\phi$ so that the assumptions in Proposition \ref{prop: KE core} are satisfied for Lebesgue measure $\mu$. We call it the Neumann Admissibility Condition (or NAC), and assuming it, Theorem \ref{theorem: NAC assures KE} calculates the top Lyapunov exponent of products of random matrices explicitly.

Let us recall the setting. Let $b \geq 2$ be an integer, and $\mu$ be a probability measure on $\mathbb{T} = \mathbb{R}/\mathbb{Z}$, invariant and ergodic for $T_b x = bx \mathrm{ \hspace{5pt} mod \hspace{3pt}} 1$. We consider $b$ number of non-zero $2 \times 2$ matrices with non-negative entries, $A_0, \ldots, A_{b-1}$. The associated M\"obius IFS $\{f_i\}_{i=0}^{b-1}$ is defined with
\begin{equation} \label{eq: definition of associated f}
f_i(x) =
\begin{pmatrix}
\alpha(i) & \beta(i) \\
\gamma(i) & \delta(i)
\end{pmatrix}(x)
:=
\begin{pmatrix}
p-q-r+s & p-q+r-s \\
p+q-r-s & p+q+r+s
\end{pmatrix}(x)
\text{\quad for \quad}
A_i =
\begin{pmatrix}
p & q \\
r & s
\end{pmatrix}.
\end{equation}

Suppose $\phi: [-1, 1] \to (-1, 1)$ is a M\"obius transformation,
\[
\phi(x) =
\begin{pmatrix}
A & B \\
C & D
\end{pmatrix}(x).
\]
We define the following constants. 
\begin{equation} \label{eq: definition of converted Mobius}
\begin{pmatrix}
p_1(i) & p_2(i) \\
q_1(i) & q_2(i)
\end{pmatrix}
=
\begin{pmatrix}
A & B \\
C & D
\end{pmatrix}
\begin{pmatrix}
\alpha(i) & \beta(i) \\
\gamma(i) & \delta(i)
\end{pmatrix}
{\begin{pmatrix}
A & B \\
C & D
\end{pmatrix}}^{-1},
\end{equation}
\begin{equation*}
u_1(i) = \left| \frac{p_2(i)}{q_2(i)} \right|, \hspace{5pt} u_2(i) = \left| \frac{p_2(i)}{q_2(i)} - \frac{p_1(i)}{q_1(i)} \right|, \hspace{5pt} u_3(i) = \left| \frac{q_1(i)}{q_2(i)} \right|, \hspace{5pt} m(i) = \max \left\{ u_1(i), \, u_2(i) \right\},
\end{equation*}
and
\begin{align*}
& \rho_1 = \max_{\substack{0 \leq i \leq b-1 \\ i: \, \mu\text{-valid}}} \left\{ u_1(i) \right\}, \hspace{5pt}
\rho_2 = \max_{\substack{0 \leq i \leq b-1 \\ i: \, \mu\text{-valid}}} \left\{ m(i) \right\}, \hspace{5pt} & \\
& \rho_3 = \max_{\substack{0 \leq i \leq b-1 \\ i: \, \mu\text{-valid}}} \left\{ \frac{m(i) \, u_3(i)}{\big( 1-m(i) \big)^2} \right\}, \hspace{5pt}
\rho_4 = \max_{\substack{0 \leq i \leq b-1 \\ i: \, \mu\text{-valid}}} \left\{ \left| \frac{u_3(i)}{1-m(i)} \right| \right\}. &
\end{align*}
We consider the following condition on $\phi$:
\begin{flalign} \label{condition: NAC}
 &\text{{\it{\textbf{Neumann Admissibility Condition (NAC)}}}: \quad}
 \rho_1 < \frac{1}{2}, \hspace{10pt}
\max \left\{
\rho_2, \, \rho_3, \, \rho_4
\right\} < 1, & \nonumber \\
& \hspace{55pt} \left| \frac{C}{A} \right| < 1, \text{\hspace{3pt} and for every $\mu$-valid index $0 \leq i \leq b-1$, \hspace{3pt}} \left| \frac{ C \delta(i) -D \gamma(i) }{A \delta(i) - B \gamma(i) } \right| < 1. &
\end{flalign}

The following theorem summarizes the result we obtained and enables us to compute the top Lyapunov exponent $\lambda$ under NAC.

\begin{theorem} \label{theorem: NAC assures KE}
Let $b \geq 2$ be an integer, and $\mu$ a probability measure on $\mathbb{T} = \mathbb{R}/\mathbb{Z}$, invariant and ergodic for $T_b x = bx \mathrm{ \hspace{5pt} mod \hspace{3pt}} 1$. Consider real matrices with non-negative entries, $A_0, \ldots, A_{b-1} \in M_2(\mathbb{R}_{\geq 0})$. Suppose there is a M\"obius transformation $\phi: [-1, 1] \to (-1, 1)$,
\[
\phi(x) = \frac{Ax+B}{Cx+D},
\]
satisfying the Neumann Admissibility Condition. Then $\{ A_0, \ldots, A_{b-1} \}$ is Kernel-expandable with respect to $\mu$, and we can define a linear operator $T: \ell^\infty(\mathbb{N}_0) \to \ell^\infty(\mathbb{N}_0)$ and $a \in \ell^\infty(\mathbb{N}_0)$ such that we have the following convergent series. $\left( \text{Here, $(\cdot)_0$ is the $0$-th entry.} \right)$
\begin{equation} \label{eq: Neumann series expression of lambda}
\lambda = \sum_{n=0}^\infty \left( T^n a \right)_0.
\end{equation}
In explicit terms, $T = (b_{k, n})_{k, n \in \mathbb{N}_0}$ is defined $\left( \text{using equation \eqref{eq: definition of converted Mobius},} \right)$ by 
\begin{align*}
&b_{k, 0} = 0  \quad \text{for all $k \in \mathbb{N}_0$,} & \nonumber \\
&b_{0, n} = \sum_{i=0}^{b-1} \, \mu(t_1 = i) \, \left( \frac{p_2(i)}{q_2(i)} \right)^n \quad \text{for $n > 0$, and} & \nonumber \\
&b_{k, n} = \sum_{i=0}^{b-1} \, \mu(t_1 = i) \, \left( \frac{p_2(i)}{q_2(i)} \right)^{n} \, \left(-\frac{q_1(i)}{q_2(i)} \right)^k \, \left\{ \sum_{\ell=1}^{\min\{k,n\}} \binom{n}{\ell}\binom{k-1}{\ell-1}
\left(1-\frac{p_1(i) \, q_2(i)}{p_2(i) \, q_1(i)}\right)^{\ell} \right\}&
\end{align*} \\[-40pt]
\begin{flalign} \label{eq: explicit definition of T}
&& \text{for $k, n>0$.}
\end{flalign}
The vector $a = \left( a_n\right)_{n \in \mathbb{N}_0} \in \ell^\infty(\mathbb{N}_0)$ is defined as $\left( \text{using equation \eqref{eq: definition of associated f},} \right)$
\[
a_0 = \sum_{i = 0}^{b-1} \, \mu(t_1 = i) \,  \log{ \left( \frac{ A \delta(i) - B \gamma(i) }{ 2A } \right) },
\]
and
\[
a_n = \frac{1}{n} \left\{ \left( \frac{C}{A} \right)^n - \sum_{i=0}^{b-1} \, \mu(t_1 = i) \, \left( \frac{ C \delta(i) -D \gamma(i) }{A \delta(i) - B \gamma(i) } \right)^n \right\}.
\]
\end{theorem}

What makes this theorem distinct is that it applies even when the $\mu$-stationary measures are continuous, and the associated algorithm computes in polylogarithmic time. The recurring method is limited in scope; the Neumann-series method applies in substantially broader situations. (For the evaluation of tails of the sum in equation \eqref{eq: Neumann series expression of lambda} and discussion on computational cost, see subsection \S\S \ref{subsection: Evaluation of Neumann series method}.) As stated in the introduction, we can indeed find a $\phi$ satisfying NAC in several examples.

The proof of this theorem will be given after two Lemmas. Recall that the infinite matrix $T = (b_{k, n})_{k, n \in \mathbb{N}_0}$ was defined in the previous subsection to satisfy
\[
\sum_{i = 0}^{b-1} \, \mu(t_1 = i) \, ( \phi \circ f_i )^n = \sum_{k=0}^\infty b_{k,n} \phi^k \quad \left( \text{ for $n > 0$ } \right).
\]
Let us confirm that the definition of $T$ in equation \eqref{eq: explicit definition of T} is coherent, and that the first assumption $(1)$ in Proposition \ref{prop: KE core} is satisfied under NAC, specifically, by the first line of equation \eqref{condition: NAC}.

\begin{lemma} \label{lemma: NAC asserts cond 1}
Under the setting of Theorem \ref{theorem: NAC assures KE}, the entries $b_{k,n}$ defined by equation \eqref{eq: explicit definition of T} satisfy the following expansion for every natural number $n$\textup{:}
\begin{equation} \label{eq: T is indeed an expansion of compositions}
\sum_{i = 0}^{b-1} \, \mu(t_1 = i) \, ( \phi \circ f_i )^n = \sum_{k=0}^\infty b_{k,n} \phi^k.
\end{equation}
Also, the infinite matrix $T = (b_{k, n})_{k, n \in \mathbb{N}_0}$ is a well-defined operator $T: \ell^\infty(\mathbb{N}_0) \to \ell^\infty(\mathbb{N}_0)$, and its operator norm $\norm{T}$ satisfies
\begin{equation} \label{eq: the bound for the operator norm}
\norm{T} \leq \max \left\{ \frac{ \rho_1 }{ 1 - \rho_1 }, \, \rho_3 \right\} < 1.
\end{equation}
\end{lemma}

\begin{proof}

The composition $\phi \circ f_i$ is again a M\"obius transformation, and can be expanded as a power series of $\phi$ as follows.
\begin{flalign*}
& ( \phi \circ f_i )^n
=
\left( {\begin{pmatrix}
A & B \\
C & D
\end{pmatrix}
\begin{pmatrix}
\alpha(i) & \beta(i) \\
\gamma(i) & \delta(i)
\end{pmatrix}
{\begin{pmatrix}
A & B \\
C & D
\end{pmatrix}}^{-1}}\big( \phi \big) \right)^n &
\end{flalign*} \vspace{-10pt}
\begin{flalign*}
&\phantom{( \phi \circ f_i )^n}
=
\left( \begin{pmatrix}
p_1(i) & p_2(i) \\
q_1(i) & q_2(i)
\end{pmatrix} \big( \phi \big) \right)^n &
\end{flalign*} \vspace{-10pt}
\begin{flalign*}
&\phantom{( \phi \circ f_i )^n}
= \Bigg( \frac{p_2(i)}{q_2(i)} + \Big( \frac{p_2(i)}{q_2(i)} - \frac{p_1(i)}{q_1(i)} \Big) \sum_{m=1}^\infty \Big(-\frac{q_1(i)}{q_2(i)}\phi\Big)^m \Bigg)^n &
\end{flalign*} \vspace{-10pt}
\begin{flalign*}
& \phantom{( \phi \circ f_i )^n}
= \left( \frac{p_2(i)}{q_2(i)} \right)^n \left\{
1 + \sum_{\ell=1}^n \binom{n}{\ell} \left(1-\frac{p_1(i) \, q_2(i)}{p_2(i) \, q_1(i)}\right)^{\ell}
\Bigg(\sum_{m=1}^\infty \left( -\frac{q_1(i)}{q_2(i)}\phi \right)^m\Bigg)^\ell \right\} &
\end{flalign*} \vspace{-10pt}
\begin{flalign*}
& \phantom{( \phi \circ f_i )^n}
= \left( \frac{p_2(i)}{q_2(i)} \right)^n \left\{
1 + \sum_{\ell=1}^n \binom{n}{\ell}
\left(1-\frac{p_1(i) \,  q_2(i)}{p_2(i) \, q_1(i)}\right)^{\ell}
\sum_{m_1,\dots,m_\ell\ge1} \prod_{j=1}^\ell \left(-\frac{q_1(i)}{q_2(i)}\phi \right)^{m_j} \right\}  &
\end{flalign*} \vspace{-10pt}
\begin{flalign} \label{equation: calculation of power of composition}
& \phantom{( \phi \circ f_i )^n}
= \left( \frac{p_2(i)}{q_2(i)} \right)^n \left\{
1 + \sum_{k=1}^\infty \sum_{\ell=1}^{\min\{k,n\}} \binom{n}{\ell}\binom{k-1}{\ell-1}
\left(1-\frac{p_1(i) \, q_2(i)}{p_2(i) \, q_1(i)}\right)^{\ell} \left(-\frac{q_1(i)}{q_2(i)}\phi \right)^k \right\}. &
\end{flalign}
The convergence is assured by NAC. This proves equation \eqref{eq: T is indeed an expansion of compositions}.

Next, let $y = (y_n)_{n=0}^\infty \in \ell^\infty(\mathbb{N}_0)$. By $\rho_1 < \frac{1}{2}$,
\begin{align*}
\left| (Ty)_0 \right| 
&= \left| \sum_{n=1}^\infty \max_{\substack{0 \leq i \leq b-1 \\ i: \, \mu\text{-valid}}}
\left\{
\left( \frac{p_2(i)}{q_2(i)} \right)^n
\right\} \, y_n \right| \leq \norm{y}_\infty \frac{\rho_1}{1- \rho_1}.
\end{align*}
For any $k \in \mathbb{N}$,
\begin{align*}
\left| (Ty)_k \right|
&= \max_{\substack{0 \leq i \leq b-1 \\ i: \, \mu\text{-valid}}} \left| \sum_{n=1}^\infty \left( \frac{p_2(i)}{q_2(i)} \right)^n \hspace{3pt} \sum_{\ell=1}^{\min\{k,n\}} \binom{n}{\ell}\binom{k-1}{\ell-1}
\left(1-\frac{p_1(i) \,  q_2(i)}{p_2(i) \, q_1(i)}\right)^{\ell}\Big(-\frac{q_1(i)}{q_2(i)} \Big)^k y_n \right| \\
&\leq \norm{y}_\infty \max_{\substack{0 \leq i \leq b-1 \\ i: \, \mu\text{-valid}}} {u_3(i)}^k \sum_{n=1}^\infty \sum_{\ell=1}^{\min\{k,n\}} \binom{n}{\ell}\binom{k-1}{\ell-1}
{u_1(i)}^\ell \, {u_2(i)}^{n-\ell} \\
&\leq \norm{y}_\infty \max_{\substack{0 \leq i \leq b-1 \\ i: \, \mu\text{-valid}}} {u_3(i)}^k \sum_{n=1}^\infty {m(i)}^n \sum_{\ell=1}^{\min\{k,n\}} \binom{n}{\ell}\binom{k-1}{\ell-1}.
\end{align*}
Here, by Vandermonde's identity,
\begin{equation*}
\sum_{\ell=1}^{\mathrm{min}\{k, \, n\}} \binom{n}{\ell} \binom{k-1}{\ell-1} = \binom{n+k-1}{k}.
\end{equation*}
Thus, by $\max \left\{
\rho_2, \rho_3, \rho_4
\right\} < 1$,
\begin{align*}
\left| (Ty)_k \right|
&\leq \norm{y}_\infty \max_{\substack{0 \leq i \leq b-1 \\ i: \, \mu\text{-valid}}} {u_3(i)}^k \sum_{n=1}^\infty \binom{n+k-1}{k} {m(i)}^n \\
&= \norm{y}_\infty \max_{\substack{0 \leq i \leq b-1 \\ i: \, \mu\text{-valid}}} \frac{m(i) \, u_3(i)}{\big( 1-m(i) \big)^2} \left( \frac{u_3(i)}{1-m(i)} \right)^{k-1}
\leq \norm{y}_\infty \, \rho_3 \, \rho_4^{\, k-1}.
\end{align*}
Therefore, $T: \ell^\infty(\mathbb{N}_0) \to \ell^\infty(\mathbb{N}_0)$ is well-defined and it satisfies the inequality \eqref{eq: the bound for the operator norm}
\end{proof}

The next lemma states that the assumption $(2)$ in Proposition \ref{prop: KE core} is satisfied by NAC, specifically, by the second line of equation \eqref{condition: NAC}.

\begin{lemma} \label{lemma: NAC asserts cond 2}
Under the setting of \hspace{1pt}Theorem \ref{theorem: NAC assures KE}, for each $\mu$-valid $0 \leq i \leq b-1$, there is $a^{(i)} = \left( a_n^{(i)} \right)_{n=0}^\infty \in \ell^\infty(\mathbb{N}_0)$ satisfying the following equality with uniformly convergent series.
\begin{equation*}
\log \norm{
\begin{pmatrix}
\frac{1+x}{2}, & \hspace{-4pt} \frac{1-x}{2}
\end{pmatrix}
\, A_i }_1
=
\sum_{n=0}^\infty a_n^{(i)} \phi(x)^n  \quad \left(\, x \in [-1, 1] \, \right).
\end{equation*}
More precisely, for each $\mu$-valid $0 \leq i \leq b-1$, the vector $a^{(i)}$ is given by
\begin{align*}
& a^{(i)}_0 = \log{ \left( \frac{ A \delta(i) - B \gamma(i) }{ 2A } \right) }
\text{, \, and for $n>0$, \, }
a^{(i)}_n = \frac{1}{n} \left\{ \left( \frac{C}{A} \right)^n -  \left( \frac{ C \delta(i) -D \gamma(i) }{A \delta(i) - B \gamma(i) } \right)^n \right\}. &
\end{align*}
\end{lemma}

\begin{proof}
For $\mu$-valid $0\leq i \leq b-1$, the integrand in \eqref{eq: lambda is the sup double int} is expressed using $f_i$ as
\begin{equation} \label{eq: integrand represented with f}
\log \norm{
\begin{pmatrix}
\frac{1+x}{2}, & \hspace{-4pt} \frac{1-x}{2}
\end{pmatrix}
\, A_i }_1
=
\log \left( \frac{\delta(i)}{2} \right) + \log \left( 1 + \frac{\gamma(i)}{\delta(i)} x \right).
\end{equation}
We can expand the second term using powers of $\phi$, as follows.
\[
\log{ \left( 1 + \frac{\gamma(i)}{\delta(i)} x \right) } = \log{ \frac{A \delta(i) - B \gamma(i) }{A \delta(i)} } + \log{ \left( 1- \frac{ C \delta(i) -D \gamma(i) }{A \delta(i) - B \gamma(i) }\phi(x) \right)} - \log{ \left( 1 - \frac{C}{A} \phi(x) \right) }.
\]
Therefore, we see that NAC is sufficient for the uniform convergence of the expansion in $\phi$, and the coefficients are as stated.
\end{proof}

Now, we are ready to prove Theorem \ref{theorem: NAC assures KE}.

\begin{proof}[Proof of Theorem \ref{theorem: NAC assures KE}] \, \\[-15pt]

By Lemma \ref{lemma: NAC asserts cond 1} and Lemma \ref{lemma: NAC asserts cond 2}, the assumptions $(1)$ and $(2)$ in Proposition \ref{prop: KE core} are satisfied. Then, the conclusions follow from Proposition \ref{prop: KE core}.
\end{proof}

\subsection{Kernel expandability through Neumann series method}
Here we will prove, using the Neumann series method introduced earlier, that some pairs of $(D_1, D_2)$ are indeed Kernel-expandable. Our goal is to prove the following theorem stated in the introduction \S\S\ref{subsection: Main results with examples}.

\begin{retheorem}{theorem: toothless twins}
Let $b \geq 7$ be a natural number and $\mu$ Lebesgue measure. Consider $\tau, u \in \{0, 1, \ldots, b-1\}$ and let
\begin{align*}
& D_1 = \{0, 1, \ldots, b-1\} \backslash \{ \tau \}, \\
& D_2 = \{0, 1, \ldots, b-1\} \backslash \{u\}.
\end{align*}
Then, the pair $(D_1, D_2)$ is either degenerate or Kernel-expandable. 
More precisely, it is Kernel-expandable if both $\tau \notin \{0, b-1\}$ and $\tau + u = b-1$. Otherwise, it is degenerate.
\end{retheorem}

Let us first deal with the Kernel-expandable cases. The strategy is to give a specific $\phi$ and to check NAC for it.

\begin{lemma} \label{lemma: Kernel expandability of toothless twins}
Take an integer $b \geq 7$, and let $\mu$ be Lebesgue measure. Let $\tau \in \{1, 2, \ldots, b-2\}$, and
\begin{align*}
& D_1 = \{0, 1, \ldots, b-1\} \backslash \{ \tau \}, \\
& D_2 = \{0, 1, \ldots, b-1\} \backslash \{b-1-\tau\}.
\end{align*}
Then, the pair $(D_1, D_2)$ is Kernel-expandable.
\end{lemma}

\begin{proof}
Suppose $\frac{b-1}{2} < \tau \leq b-2$. (The proof for when $1 \leq \tau \leq \frac{b-1}{2}$ works in the same way.) Then,
\begin{equation*}
\# \big( \left( D_1 + s \right) \cap D_2 \big) = \left\{
\begin{array}{ll}
b-s-2
& \text{\quad for \, $0 \le s < b-\tau$,} \\
b-s
& \text{\quad for \, $b-\tau \le s \le b$.}
\end{array}
\right.
\end{equation*}

\begin{equation*}
\# \big( \left( D_1 + s \right) \cap \left(D_2 + b \right) \big) = \left\{
\begin{array}{ll}
s
& \text{\quad for \, $0 \le s < b-\tau$,} \\
s-2
& \text{\quad for \, $b-\tau \le s < 2b-2\tau-1$,} \\
s-1
& \text{\quad for \, $s = 2b-2\tau-1$,} \\
s-2
& \text{\quad for \, $2b-2\tau-1 < s \le b$.}
\end{array}
\right.
\end{equation*}
Therefore, the matrices $A_i$ are calculated as follows.
\begin{equation*}
A_i =
\begin{cases}
\vcenter{\hbox{$\begin{pmatrix}
\colA{b-i-2} & \colB{i}\\
\colA{b-i-3} & \colB{i+1}
\end{pmatrix}$}}
& \text{for \,} 0 \le i < b-\tau-1,\\[16pt]

\vcenter{\hbox{$\begin{pmatrix}
\colA{\tau-1} & \colB{b-\tau-1}\\
\colA{\tau}   & \colB{b-\tau-2}
\end{pmatrix}$}}
& \text{for \,} i = b-\tau-1,\\[16pt]

\vcenter{\hbox{$\begin{pmatrix}
\colA{b-i}   & \colB{i-2}\\
\colA{b-i-1} & \colB{i-1}
\end{pmatrix}$}}
& \text{for \,} b-\tau \le i < 2b-2\tau-2,\\[16pt]

\vcenter{\hbox{$\begin{pmatrix}
\colA{2\tau-b+2} & \colB{2b-2\tau-4}\\
\colA{2\tau-b+1} & \colB{2b-2\tau-2}
\end{pmatrix}$}}
& \text{for \,} i = 2b-2\tau-2,\\[16pt]

\vcenter{\hbox{$\begin{pmatrix}
\colA{2\tau-b+1} & \colB{2b-2\tau-2}\\
\colA{2\tau-b}   & \colB{2b-2\tau-2}
\end{pmatrix}$}}
& \text{for \,} i = 2b-2\tau-1,\\[16pt]

\vcenter{\hbox{$\begin{pmatrix}
\colA{b-i}   & \colB{i-2}\\
\colA{b-i-1} & \colB{i-1}
\end{pmatrix}$}}
& \text{for \,} 2b-2\tau-1 < i \le b-1.
\end{cases}
\end{equation*}
In turn, using equation \eqref{eq: induced map by Ai}, associated maps $f_i$ are
\begin{equation} \label{eq: calculated f}
f_i =
\begin{cases}
\vcenter{\hbox{$\begin{pmatrix}
\colA{2} & \colB{2b-4i-6 \hspace{4pt}}\\
\colA{0} & \colB{2b-4 \hspace{4pt}}
\end{pmatrix}$}}
& \text{for \,} 0 \le i < b-\tau-1,\\[16pt]

\vcenter{\hbox{$\begin{pmatrix}
\colA{-2} & \colB{-2b+4\tau+2 \hspace{4pt}}\\
\colA{0}  & \colB{2b-4 \hspace{4pt}}
\end{pmatrix}$}}
& \text{for \,} i = b-\tau-1,\\[16pt]

\vcenter{\hbox{$\begin{pmatrix}
\colA{2} & \colB{2b-4i+2 \hspace{4pt}}\\
\colA{0} & \colB{2b-4 \hspace{4pt}}
\end{pmatrix}$}}
& \text{for \,} b-\tau \le i < 2b-2\tau-2,\\[16pt]

\vcenter{\hbox{$\begin{pmatrix}
\colA{3} & \colB{-6b+8\tau+9 \hspace{4pt}}\\
\colA{-1} & \colB{2b-3 \hspace{4pt}}
\end{pmatrix}$}}
& \text{for \,} i = 2b-2\tau-2,\\[16pt]

\vcenter{\hbox{$\begin{pmatrix}
\colA{1} & \colB{-6b+8\tau+5 \hspace{4pt}}\\
\colA{1} & \colB{2b-3 \hspace{4pt}}
\end{pmatrix}$}}
& \text{for \,} i = 2b-2\tau-1,\\[16pt]

\vcenter{\hbox{$\begin{pmatrix}
\colA{2} & \colB{2b-4i+2 \hspace{4pt}}\\
\colA{0} & \colB{2b-4 \hspace{4pt}}
\end{pmatrix}$}}
& \text{for \,} 2b-2\tau-1 < i \le b-1.
\end{cases}
\end{equation}

Now, we need to find values of $A,B,C,D$ so that NAC is satisfied. This requires numerical assistance; the scripts used in this proof for computing the rigorous guaranteed bound, as well as their output, are archived at the following repository: \\
\url{https://github.com/NimaAlibabaei/On-the-intersection-of-Cantor-sets}, \\
commit: \texttt{de92aa09b1a0d9b6b61d4d7bc2fa2148afcfac8f}.

For $b \geq 22$ we choose $\phi(x) = \frac{15x+4}{5x+40}$. The second line of NAC is easily checked: under our setting, the $\gamma(i)$ in equation \eqref{eq: definition of associated f} is non-zero only for $i = 2b-2\tau-2, 2b-2\tau-1$. Thus, the second line of NAC is equivalent to
\begin{equation*}
\left| \frac{ 10b-15 \pm 40 }{ 30b-45 \pm 4 } \right| < 1 \text{\quad and \quad} \left| \frac{1}{2} \right| < 1.
\end{equation*}
The first inequality is satisfied for $b \geq 22$.

The first line of NAC is more complicated to check. We first consider the case $0 \leq i < b-\tau-1$ and change the variables with $u(i) = \frac{i}{b}$ and $v = \frac{1}{b}$. Then, $0 \leq u \leq \frac{1}{2}$, $0 \leq v \leq \frac{1}{22}$, and
\begin{align*}
\begin{pmatrix}
p_1(i) & p_2(i) \\
q_1(i) & q_2(i)
\end{pmatrix}
= \frac{80}{29}
\begin{pmatrix}
150u(i) + 865v - 95 & -450u(i) - 855v + 285 \\
50u(i) + 675v - 225 & -150u(i) -1445v + 675
\end{pmatrix}.
\end{align*}
We rigorously verified NAC for $(u, v) \in [0, \frac{1}{2}] \times [0, \frac{1}{22}]$ using validated numerical computation. In particular, we obtained
\[ \rho_1 \leq 0.43, \, \, \rho_2 \leq 0.43, \, \, \rho_3 \leq 0.44, \, \, \rho_4 \leq 0.59. \]
Hence, the desired inequalities are satisfied. Similar computations are carried for each six cases of $0 \leq i \leq b-1$ in \eqref{eq: calculated f} after appropriate change of variables.

We are left to check the case $7 \leq b \leq 21$.For these $b$ and $1 \leq \tau \leq b-2$, we rigorously proven the existence of the values of $A, B, C$ satisfying NAC for each pair of $(D_1, D_2)$, thus completing the proof.

\end{proof}

Finally, let us quickly go over the degenerate cases.

\begin{lemma} \label{lemma: extreme digits and degeneracy}
Let $\mu$ be Lebesgue measure. The pair $(D_1, D_2)$ is degenerate if one of $D_1$, $D_2$ misses one of the extreme digits $0, b-1$.
\end{lemma}

\begin{proof}
If $0 \notin D_1$, then
\begin{align*}
| (D_1 + b-1) \cap D_2 | = | (D_1 + b) \cap D_2 | = 0.
\end{align*}
Thus, with some $n,m \in \mathbb{N}_0$,
\begin{equation*}
A_{b-1} = 
\begin{pmatrix}
0 & n \\
0 & m
\end{pmatrix},
\end{equation*}
which implies that $(D_1, D_2)$ is degenerate. If $b-1 \notin D_1$,
\begin{align*}
| D_1 \cap (D_2 + b) | = | (D_1 + 1) \cap (D_2 + b) | = 0.
\end{align*}
We have, with some $n, m \in \mathbb{N}_0$,
\begin{equation*}
A_0 = 
\begin{pmatrix}
n & 0 \\
m & 0
\end{pmatrix}.
\end{equation*}
Thus $(D_1, D_2)$ is degenerate. Similar arguments work for when $D_2$ misses these digits.

\end{proof}

\begin{lemma} \label{lemma: degeneracy of toothless twins}
Let $\mu$ be Lebesgue measure. Let $\tau, u \in \{0, 1, \ldots, b-1\}$, and
\begin{align*}
& D_1 = \{0, 1, \ldots, b-1\} \backslash \{ \tau \}, \\
& D_2 = \{0, 1, \ldots, b-1\} \backslash \{u\}.
\end{align*}
Then, the pair $(D_1, D_2)$ is degenerate if one of $\{ \tau, u\}$ is an extreme digit $\{0, b-1\}$ or $\tau + u \ne b-1$.
\end{lemma}

\begin{proof}
By Lemma \ref{lemma: extreme digits and degeneracy}, if $\tau$ or $u \in \{0, b-1\}$, the pair $(D_1, D_2)$ is degenerate. Suppose $\tau, u \notin \{0, b-1\}$ and $\tau + u \ne b-1$. We will see that $A_u$ is always degenerate.

Suppose either (i) $\tau < u$ and $b - 1 - u < \tau$, 
(ii) $u < \tau$ and $b - 1 - \tau < u$, or 
(iii) $\tau = u$ and $b - 1 < 2u$. 
In any of the above case, it can be checked that
\begin{align*}
&| (D_1 + u) \cap D_2 | = b-u-1, \quad | (D_1 + u) \cap (D_2 + b) | = u-1 \\
&| (D_1 + u+1) \cap D_2 | = b-u-1, \quad | (D_1 + u+1) \cap (D_2 + b) | = u-1.
\end{align*}
Thus,
\begin{equation*}
A_u = 
\begin{pmatrix}
b-u-1 & u-1 \\
b-u-1 & u-1
\end{pmatrix}.
\end{equation*}
Similarly, when either (iv) $\tau < u$ and $\tau < b - 1 - u$, 
(v) $u < \tau$ and $u < b - 1 - \tau$, or 
(vi) $\tau = u$ and $2u < b - 1$,
we have
\begin{equation*}
A_u = 
\begin{pmatrix}
b-u-2 & u \\
b-u-2 & u
\end{pmatrix}.
\end{equation*}
Therefore, the pair $(D_1, D_2)$ is degenerate under this assumption.
\end{proof}

Combining the Lemma \ref{lemma: Kernel expandability of toothless twins} and \ref{lemma: degeneracy of toothless twins} completes the proof of Theorem \ref{theorem: toothless twins}.

\subsection{Evaluation of Neumann series method} \label{subsection: Evaluation of Neumann series method}
We now know that the Neumann series method can indeed be used to establish Kernel-expandability. To actually compute the dimension using the expression obtained in Theorem~\ref{theorem: NAC assures KE}, we must quantify the contribution of the truncated tails.

In this section we always assume NAC in equation \eqref{condition: NAC}. Under the notations of NAC, let
\begin{equation} \label{eq: definition of max r}
r = \max \left\{ \max_{\substack{0 \leq i \leq b-1 \\ i: \, \mu\text{-valid}}} \left\{ \left| \frac{ C \delta(i) -D \gamma(i) }{A \delta(i) - B \gamma(i) } \right| \right\}, \hspace{3pt} \left| \frac{C}{A} \right|, \hspace{3pt} \rho_4 \right\}.
\end{equation}
Then $0 < r < 1$. By Lemma \ref{lemma: NAC asserts cond 2}, in terms of evaluating $\sum T^n y$, we may only consider \ $y = (y_n)_{n=0}^\infty \in \ell^\infty(\mathbb{N}_0)$ with $\abs{y_n} \leq r^{n-1}$ for all $n \in \mathbb{N}$. Since constants penetrate the integral, we can always assume $y_0 = 0$.

Let
\[
E = \max \left\{ \frac{1}{1- r \rho_1}, \rho_3 \right\}, \quad R_N = E r^{N-1}.
\]
We have $\lim_{N \to \infty} R_N = 0$.
 
For an infinite matrix $T = (b_{k, n})_{k, n \in \mathbb{N}_0}$, define $T_N = (c_{k,n})_{k, n \in \mathbb{N}_0}$ by $c_{k,n} = b_{k,n}$ when $k, n \leq N-1$ and $0$ otherwise.

\begin{lemma} \label{lemma: iterated evaluation}
Under NAC, for any natural number $n$,
\begin{align*}
\left\| T^n y - {T_N}^n y \right\|_\infty \leq n \, R_N.
\end{align*}
\end{lemma}

\begin{proof}
First, consider $n=1$.
\begin{equation*}
\left| \left( Ty - T_N y \right)_0 \right|
\leq \max_{\substack{0 \leq i \leq b-1 \\ i: \, \mu\text{-valid}}} \frac{ \big( r \, u_1(i)\big)^{N-1}}{1- r \, u_1(i)}
\leq \frac{ \left(r \, \rho_1\right)^{N-1} }{1- r \, \rho_1} \leq R_N,
\end{equation*}
and for $0 < k < N$,
\begin{align*}
\left| \left( Ty - T_N y \right)_k \right|
&\leq \max_{\substack{0 \leq i \leq b-1 \\ i: \, \mu\text{-valid}}} \, {u_3(i)}^k \sum_{n=N}^\infty \binom{n+k-1}{k} m(i)^n r^{n-1} \\
&\leq \max_{\substack{0 \leq i \leq b-1 \\ i: \, \mu\text{-valid}}} \, r^{N-1} \, {u_3(i)}^k \sum_{n=1}^\infty \binom{n+k-1}{k} m(i)^n \\
&\leq \max_{\substack{0 \leq i \leq b-1 \\ i: \, \mu\text{-valid}}} \frac{ \, m(i) \, u_3(i)}{\big(1-m(i) \big)^2} \, r^{N-1} \left( \frac{u_3(i)}{1-m(i)} \right)^{k-1} \\
&\leq \, \rho_3 \, r^{N-1} \, \rho_4^{k-1} \leq R_N.
\end{align*}
Finally, for $k \geq N$,
\begin{align}
\left| \left( Ty - T_N y \right)_k \right|
= \left| \left( Ty \right)_k \right|
&\leq \max_{\substack{0 \leq i \leq b-1 \\ i: \, \mu\text{-valid}}} \, {u_3(i)}^k \sum_{n=1}^\infty \binom{n+k-1}{k} m(i)^n r^{n-1}  \nonumber \\
&\leq \max_{\substack{0 \leq i \leq b-1 \\ i: \, \mu\text{-valid}}} \frac{\, m(i) \, u_3(i)}{\big( 1-m(i) \big)^2} \left( \frac{u_3(i)}{1-m(i)} \right)^{N-1} \nonumber \\
&\leq \, \rho_3 \, {\rho_4}^{N-1} \leq R_N \label{eq: evaluation of k-th Ty}.
\end{align}
Therefore, \vspace{-15pt}
\begin{align*}
\left\| Ty - T_N y \right\|_\infty \leq R_N.
\end{align*}

We notice that the sequence $T_N y \in \ell^\infty(\mathbb{N}_0)$ is also bounded by $(r^n)_n$; as in inequality \eqref{eq: evaluation of k-th Ty},
\begin{equation*}
\left| ( T_N y)_k \right|
\leq \max_{\substack{0 \leq i \leq b-1 \\ i: \, \mu\text{-valid}}} \frac{ \, m(i) \, u_3(i)}{\big( 1-m(i) \big)^2} \left( \frac{u_3(i)}{1-m(i)} \right)^{k-1}
\leq  \, \rho_4^{\, k-1} \leq r^{k-1},
\end{equation*}
Note that $(T_N y)_0$ might be positive, but since the left-most column of $T$ is $0$ everywhere, this does not carry on  to $T_N(T_N y)$. Thus, for any $n, k > 0$,
\begin{equation*}
\left( T_N \left( {T_N}^n y \right) \right)_k \leq \, r^{\, k-1},
\end{equation*}
Using the initial argument again, we have for any $n > 1$,
\begin{equation*}
\left\| T \left( {T_N}^n y \right) - T_N \left( {T_N}^n y \right) \right\|_\infty \leq R_N,
\end{equation*}

The desired inequality follows by induction. For any natural number $n > 2$,
\begin{align*}
\left\| T^n y - {T_N}^n y \right\|_\infty
&\leq \left\| T^n y - T \left( {T_N}^{n-1} y \right) \right\|_\infty + \left\| T \left( {T_N}^{n-1} y \right) - T_N \left( {T_N}^{n-1} y \right) \right\|_\infty \\
&\leq \norm{T} (n-1) R_N + R_N \\
&\leq n \, R_N.
\end{align*}

\end{proof}

\begin{proposition} \label{proposition: evaluation of Neumann series}
Assume NAC. Let $r$ be defined by equation \eqref{eq: definition of max r}. Let $y = (y_n)_{n=0}^\infty \in \ell^\infty(\mathbb{N}_0)$ with $\abs{y_n} \leq  r^{n-1}$ for all natural numbers $n > 0$, and $y_0 = 0$. Then, for a given $\vep > 0$, there are natural numbers $M$ and $N$ such that
\begin{equation*}
\left\| \sum_{n=0}^\infty T^n y - \sum_{n=0}^M {T_N}^n y \right\|_\infty \leq \vep.
\end{equation*}
\end{proposition}

\begin{proof}
For a natural number $M$, we have
\begin{equation*}
\left\| \sum_{n = M+1}^\infty T^n y \right\|_\infty \leq \sum_{n=M+1}^\infty \norm{T}^n \leq \frac{\norm{T}^{M+1} }{ 1 - \norm{T} }.
\end{equation*}
Then, take $M$ so that
\begin{equation*}
\frac{ \norm{T}^{M+1} }{ 1 - \norm{T} } < \frac{\vep}{2}.
\end{equation*}
Also, by Lemma \ref{lemma: iterated evaluation},
\begin{equation*}
\left\| \sum_{n=0}^M T^n y - \sum_{n=0}^M {T_N}^n y \right\|_\infty \leq \frac{ M(M+1) }{2} R_N.
\end{equation*}
Thus, take $N$ so that
\begin{equation*}
\frac{ M(M+1) }{2} R_N < \frac{\vep}{2}.
\end{equation*}

\end{proof}

\begin{remark} \label{remark: computational cost}
From the inequalities in the proof, one has the asymptotic scaling
\[
M = O \! \left( \log{ ( 1/\vep ) } \right), \quad N = O \! \left( \log{ ( 1/\vep ) } \right).
\]
The cost to build the truncated matrix $T_N$ by using equation \eqref{eq: explicit definition of T} is
\[
O \! \left( \sum_{k=1}^N \sum_{n=1}^N \min{(k,n)}  \right)
=
O \! \left( N^3  \right).
\]
In addition, the cost of forming $T_N^m y$ for all $0 \leq m \leq M$ by repeated matrix–vector multiplication is $O \! \left( M N^2 \right)$ arithmetic operations; hence, one obtains the polylogarithmic-time cost of
\[
O \! \left( N^3  \right) + O \! \left( M N^2 \right)
=
O \! \left( \Big( \log{ \left( 1/\vep \right) } \Big)^3  \right).
\]
\end{remark}

\section*{Acknowledgement}
This paper grew out of countless conversations with my advisor, Masaki Tsukamoto — his persistence and high standards guided me at every step.

I also want to thank the friends and family who kept me grounded during the writing process. I am grateful to the scholarly community whose generosity with time and support made this project possible.

I had help from ChatGPT (OpenAI) to develop the codes used in the proof. All codes were reviewed and integrated by me, who is responsible for the final implementation and results.

\vspace{0.5cm}

\address{
Department of Mathematics, Kyoto University, Kyoto 606-8501, Japan}

\textit{E-mail}: \texttt{alibabaei.nima.28c@st.kyoto-u.ac.jp}

This work was supported by JSPS KAKENHI Grant Number 25KJ1473.

\end{document}